\title{Nearly optimal embedding of trees}

\author
{Benny Sudakov\thanks{Department of Mathematics, Princeton,
Princeton, NJ 08544, USA. Email: {\tt bsudakov@math.princeton.edu}. Research
supported in part by NSF CAREER award DMS-0546523, NSF grant
DMS-0355497 and by a USA-Israeli BSF grant.}
\and Jan Vondr\'ak \thanks{
Department of Mathematics, Princeton,
Princeton, NJ 08544, USA.
Email: {\tt jvondrak@math.princeton.edu}.}
}

\documentclass[11pt]{article}
\usepackage{amsfonts}
\usepackage{amsmath}
\usepackage{latexsym}
\usepackage{epsfig}
\usepackage{graphicx}
\usepackage{amssymb}

\newtheorem{theorem}{Theorem}[section]
\newtheorem{lemma}[theorem]{Lemma}
\newtheorem{proposition}[theorem]{Proposition}

\newtheorem{definition}[theorem]{Definition}

\def\Pr{{\mathbb P}}
\def\E{{\mathbb E}}

\date{}

\begin{document}
\maketitle

\begin{abstract}
In this paper we show how to find nearly optimal
embeddings of large trees in several natural classes of graphs.
The size of the tree $T$ can be as large as a constant fraction of the size
of the graph $G$, and the maximum degree of $T$ can be close to the minimum
degree of $G$. For example, we prove that any graph of minimum degree $d$
without $4$-cycles contains every tree of size $\epsilon d^2$ and maximum
degree at most $d - 2 \epsilon d - 2$. As there exist $d$-regular graphs
without $4$-cycles of size $O(d^2)$, this result is optimal up
to constant factors. We prove similar nearly tight results for graphs
of given girth, graphs with no complete bipartite subgraph $K_{s,t}$,
random and certain pseudorandom graphs. These results are obtained using a simple 
and very natural randomized embedding algorithm, 
which can be viewed as a "self-avoiding tree-indexed random walk".

\end{abstract}

\section{Introduction}

We consider the problem of embedding a tree $T$ in a given graph $G$.
Formally, we look for an injective map $f:V(T) \rightarrow V(G)$ which
preserves the edges. We do not require that non-edges are mapped to non-edges,
i.e. the copy of $T$ in $G$ need not be induced. Our goal is to find
sufficient conditions on $G$ in order to contain all trees of certain size,
with maximum degree as large as a constant fraction (possibly approaching 1) of the 
minimum degree of $G$.

\subsection{Brief history}

The problem of embedding paths and trees in graphs has long been one
of the fundamental questions in combinatorics. This problem has been extensively studied in extremal combinatorics, in the theory of random graphs,
in connection with properties of expanders and with applications to Computer Science. 
The goal always has been to find a suitable property of a graph $G$ which guarantees that it contains all possible trees with given parameters. We describe next several examples which we think are representative and give a good overview of 
previous research in this area.

\paragraph{Extremal questions.} The basic extremal question about trees is to determine the number of edges that a graph needs to have in order to contain
all trees of given size.
It is an old folklore result that a graph $G$ of minimal degree $d$
contains every tree $T$ with $d$ edges. This can be achieved simply by embedding vertices of $T$ greedily one by one. 
Since at most $d$ vertices of $G$ are occupied at any point,
there is always enough room to embed another vertex of the tree.

An old conjecture of Erd\H{o}s and S\'os says that {\em average degree} $d$
is already sufficient to guarantee the same property.
More precisely, any graph with more than $(d-1)n/2$ edges
contains all trees with $d$ edges. A clique of size $d$ is an obvious
tight example for this conjecture. The conjecture has been proved in several 
special cases, e.g. Brandt and Dobson \cite{BD96} establish it for graphs
of girth at least $5$ ({\em girth} is the length of the shortest cycle in
a graph). In fact, they prove a stronger statement, that any such graph
of minimum degree $d/2$ and maximum degree $\Delta$ contains all trees
with $d$ edges and maximum degree at most $\Delta$.
More generally, improving an earlier result of \L uczak and Haxell \cite{HL00},
Jiang proved that any graph of girth $2k+1$ and minimum degree $d/k$ contains
all trees with $d$ edges and maximum degree at most $d/k$ \cite{Jiang01}.
For general graphs, it has been announced by Ajtai, Koml\'os, Simonovits
and Szemer\'edi \cite{AKSS} that they proved Erd\H{o}s-S\'os conjecture
for all sufficiently large trees.

A related statement, known as Loebl's $(\frac{n}{2}-\frac{n}{2}-\frac{n}{2})$
conjecture \cite{ELS94}, is that any graph on $n$ vertices,
with at least $n/2$ vertices of degree at least $n/2$,
contains all trees with at most $n/2$ edges.
Progress on this conjecture has been recently made by Yi Zhao \cite{Zhao07}.
Note that in the results discussed so far, the size of the
tree is of the same order as degrees in the graph $G$. Without assuming any
additional properties of $G$, this seems to be a natural barrier. 

\paragraph{Expanding graphs.}
Embedding trees of size much larger than the average degree of the graph
is possible in graphs satisfying certain expansion properties.
The first such result was established by P\'osa using his celebrated rotation/extension technique. 
Given a subset of vertices $X$ of a graph $G$ let
$N(X)$ denote the set of all neighbors of vertices of $X$ in $G$.
P\'osa \cite{Posa76} proved that if
$|N(X) \setminus X| \geq 2|X|-1$ for every subset $X$ of $G$ with at most $t$ vertices, then
$G$ contains a path of length $3t-2$.
This technique was extended to trees by Friedman and Pippenger \cite{FP87}.
They proved that if $|N(X)| \geq (d+1) |X|$ for all subsets of size at most 
$2t-2$, then $G$ contains every tree of size $t$ and maximum degree
at most $d$. The power of this technique is that while $T$ can have degrees
close to the minimum degree of $G$, it can be of size much larger than $d$,
depending on the expansion guarantee. On the other hand, note that these techniques cannot embed
trees of size larger than $|G|/ d$, due to the nature of the expansion property.
The result of Friedman and Pippenger has several interesting applications. 
For example, it can be used to show that for a fixed $\delta>0$, $d$ and
every $n$ there is a graph $G$ with $O(n)$ edges that, even after deletion of all but 
$\delta |E(G)|$ edges, continues to contain every tree with $n$ vertices and maximum
degree at most $d$.
This has immediate corollaries in Ramsey Theory.

The technique from \cite{FP87} also has an application for infinite graphs.
For an infinite graph $G$, its {\em Cheeger constant} is
$h(G)=\inf_X\frac{|N(X) \setminus X|}{|X|}$,
 where $X$ is a nonempty finite subset of vertices of $G$.
Using the ideas of Friedman and Pippenger, one can show
(see \cite{BS97}) that any infinite graph $G$ with Cheeger constant
$d \geq 3$ contains an infinite tree $T$ with Cheeger constant
$d-2$. Benjamini and Schramm \cite{BS97} prove a stronger result that
any infinite graph with $h(G) > 0$ contains an infinite tree with
positive Cheeger constant. They use the notion of {\em tree-indexed
random walks} to find such a tree. We will allude to this notion again later.

\paragraph{Random and pseudorandom graphs.}
The random graph $G_{n,p}$ is a probability space whose
points are graphs on a fixed set of $n$ vertices, where each pair
of vertices forms an edge, randomly and independently, with
probability $p$. For random graphs, Erd\H{o}s conjectured that with high probability,
$G_{n,d/n}$ for a fixed $d$ contains a very long path, i.e., a path of length
$(1-\alpha(d)) n$ such that $\lim_{d \rightarrow \infty} \alpha(d) = 0$.
This conjecture was proved by Ajtai, Koml\'os and Szemer\'edi \cite{AKS} and,
in a slightly weaker form, by Fernandez de la Vega \cite{FD}.
Embedding trees, however, is considerably harder. Fernandez de la Vega
\cite{Vega88} showed that there are (large) constants $a_1, a_2$ such that
$G_{n,d/n}$ contains any {\em fixed} tree $T$ of size $n / a_1$ and maximum
degree $\Delta \leq d / a_2$ w.h.p (i.e., with probability tending to 1 when $n \rightarrow \infty$). 
Note that this is much weaker than
containing all trees simultaneously, because a random graph can contain every
fixed tree w.h.p, and still miss at least one tree w.h.p.
Until recently, there was no result known on embedding all trees
simultaneously. Alon, Krivelevich and Sudakov proved in \cite{AKS}
that for any $\epsilon>0$, $G_{n, d/n}$ contains all trees of size $(1-\epsilon) n$
and maximum degree $\Delta$ such that
$$ d \geq \frac{10^6}{\epsilon} \Delta^3 \log \Delta \log^2 (2/\epsilon).$$
(All logarithms here and in the rest of this paper have natural base.)
This result is nearly tight in terms of the size of $T$, and holds
for all trees simultaneously. But it is achieved at the price of requiring that degrees in $G$
are much larger than degrees in the tree. A similar result for 
pseudorandom graphs was also proved in \cite{AKS}.
A graph $G$ is called an $(n,d,\lambda)$-graph if $G$ has $n$ vertices,
is $d$-regular (hence the largest eigenvalue of the adjacency matrix is $d$)
and the second largest eigenvalue is $\lambda$. Such graphs are known to have
good expansion and other random-like properties. Alon, Krivelevich and Sudakov proved
that any $(n,d,\lambda)$-graph such that
$$ \frac{d}{\lambda} \geq \frac{160}{\epsilon} \Delta^{5/2} \log (2/\epsilon) $$
contains all trees of size $(1-\epsilon) n$ and degrees bounded by $\Delta$.
Note that using the expansion properties of $(n,d,\lambda)$-graphs, one could
have used Friedman-Pippenger as well; however, one would not be able to embed
trees larger than $n / \Delta$ in this way.

\paragraph{Universal graphs.}
In a more general context, graphs containing all trees with given parameters
can be seen as instances of {\em universal graphs}. For a family of graphs
$\cal F$, a graph $G$ is called $\cal F$-universal, if it contains every
member of $\cal F$ as a subgraph. The construction of $\cal F$-universal
graphs for various families of subgraphs is important in applied areas
such as VLSI design, data representation and parallel computing. For trees,
a construction is known of a graph $G$ on $n$ vertices which contains
all trees with $n$ vertices and degrees bounded by $d$, such that the maximum
degree in $G$ is a function of $d$ only \cite{BCLR89}.

\subsection {Our results}

We prove several results concerning embedding trees in graphs with no short cycles,
graphs without a given complete bipartite subgraph, random graphs
and also graphs satisfying a certain pseudorandomness property.
We embed trees with parameters very close to trivial upper bounds
that cannot be exceeded: maximum degree close to the minimum
degree of $G$, and size a constant fraction of the order of
$G$~\footnote{By the {\em order} of a graph, we mean the number of vertices.
By {\em size}, we mean the number of edges. For trees, the two quantities
differ only by $1$.}
(or more precisely the minimum possible order of $G$ under given conditions).
A summary of our main results follows. Here we assume that $d$ and $n$ are sufficiently large.

\begin{enumerate}
\item For any constant $k \geq 2$, $\epsilon \leq \frac{1}{2k}$
and any graph $G$ of girth at least $2k+1$ and minimum degree $d$,
$G$ contains every tree $T$ of size $|T| \leq \frac14 \epsilon d^k$
and maximum degree $\Delta \leq (1-2\epsilon) d-2$.
\item For any $G$ of minimum degree $d$, not containing $K_{s,t}$
 (a complete bipartite graph with parts of size $s \geq t \geq 2$), 
$G$ contains every tree $T$ of size $|T| \leq \frac{1}{64 s^{1/(t-1)}}
 d^{1+\frac{1}{t-1}}$ and maximum degree $\Delta \leq \frac{1}{256} d$.
\item For a random graph $G_{n,p}$ with $d = pn \geq n^{1/k}$ for
some constant $k$, with high probability $G_{n,p}$ contains all trees
of size $O(n/k)$ and maximum degree $O(d/k)$.
\end{enumerate}

It is easy to see that any graph of girth $2k+1$ and minimum degree $d$
has at least $\Omega(d^k)$ vertices. It is a major open question to determine the smallest possible order of such graph.
For values of $k = 2,3,5$ there are known constructions obtained by 
Erd\H{o}s and R\'enyi \cite{ErdR} and Benson \cite{B} of graphs of girth $2k+1$, minimum degree $d$ and order $O(d^k)$. It is 
also widely believed that such constructions should be possible for all fixed $k$.
This implies that our first statement is tight up to constant factors for $k = 2,3,5$ and 
probably for all remaining $k$. Similarly, it is conjectured that for $s \geq t$ there are  $K_{s,t}$-free graphs  
with minimum degree $d$ which have $O(d^{1+\frac{1}{t-1}})$ vertices.
For $s > (t-1)!$,
such a construction was obtained by Alon, R\'onyai and Szabo \cite{ARS}
(modifying the construction in \cite{KRS}).
Hence, the size of the trees we are embedding in our second result is tight up to constant factors as well.
Finally, since the minimum degree of the random graph 
$G_{n,p}$ is roughly $pn$, it is easy to see that for constant $\alpha>0$
and $p=n^{-\alpha}$ we are embedding trees whose size and maximum degree is
proportional to the order and the minimum degree of $G_{n,p}$.
Thus our third result is also nearly optimal.

\subsection{Discussion}

\paragraph{Local expansion.}
Using well known results from Extremal Graph
theory, one can show that if graph $G$ contains no subgraphs isomorphic to a fixed
bipartite graph $H$ (e.g., $C_{2k}$ or $K_{s,t}$) then it has certain expansion
properties. More precisely, all small subsets of $G$ have a large boundary.
For example, if $G$ is a $C_4$-free graph with minimum degree $d$ then
all subsets of $G$ of size at most $d$ 
expand by a factor of $\Theta(d)$. Otherwise we would get a $4$-cycle
by counting the number of edges between  $S$ and its boundary $N(S) \setminus S$.
This simple observation appears to be a powerful tool in attacking
various extremal problems and was used in \cite{SV} and \cite{KS}
to resolve several conjectures about cycle lengths and clique-minors in $H$-free graphs.

Therefore, it is natural to ask whether the expansion of $H$-free graphs combined
with the result of Friedman and Pippenger can be used to embed large trees.
Recall that to embed a 
tree of size $t$ of maximum degree $d$, Friedman and Pippenger require
that sets of size up to $2t-2$ expand at least $d+1$ times.
For example, plugging this into
the observation we made on the expansion of $C_4$-free graphs only gives embedding of
trees of order $O(d)$ in such graphs. This is quite far from the bound $O(d^2)$
which can be achieved using our approach. Similarly, in graphs of girth $2k+1$,
we can embed trees of size $O(d^k)$, rather than $O(d^{k-1})$ as can be guaranteed
by using Friedman-Pippenger. Therefore, our work can be seen as an extension of the embedding results for locally expanding graphs. It shows that using structural information about $G$, rather then just local expansion, one can embed in $G$
trees of much larger size.

\paragraph{Extremal results.}
Our work sheds some light on why the Erd\H{o}s-S\'os conjecture, which we already discussed in the beginning of the introduction,  
becomes easier for graphs with no short cycles. This scenario was considered, e.g., in 
\cite{BD96, HL00, Jiang01}. In particular, assuming that graph $G$ has girth $2k+1, k\geq 2$
and minimum degree $d$, Jiang \cite{Jiang01} showed how to embed in $G$ all
trees of size $kd$ with degrees bounded by $d$. Although this is best possible,
our result implies that this 
statement can be tight only for a relatively few very special trees,
i.e., those that contain several large stars of degree $d$
or extremely close to $d$. Indeed, if we relax the degree assumption and consider trees with the maximum degree 
at most $(1-\epsilon) d$, then it is possible to embed trees of size $O(d^k)$ rather than $O(d)$.
Moreover, a careful analysis of our proof shows that it still works for
$\epsilon$ which have order of magnitude $k\frac{\log d}{d}$. Therefore even if we allow the degree of the tree to be as large as $d-ck\log d$ for some constant $c$,
we are still able to embed all trees of size $\Omega(kd^{k-1}\log d) \gg kd$.

\paragraph{Random graphs.}
It is quite easy to prove an analog of the result of Fernandez de la Vega \cite{Vega88} on the embedding of a fixed tree of size proportional to $n$
and maximum degree $O(pn)$ in the dense random graph $G_{n,p}$.
Indeed for constant $\alpha<1$ and edge probability $p = n^{-\alpha}$, this 
can be done greedily, vertex by vertex, generating the random graph simultaneously
with the embedding. On the other hand, this simple approach cannot be used to embed
all such trees with high probability, since there are too many trees to use
the union bound.
We provide the first result for simultaneous embedding of all trees of size $\Theta(n)$ and maximum degree $\Theta(pn)$,
in the random graph $G_{n,p}$ for $p = n^{-\alpha}$ and constant $\alpha<1$.
It is also interesting to compare our result with the work of Alon, 
Krivelevich, Sudakov \cite{AKS}. They embed nearly spanning trees
but with degree which is only a small power (roughly $1/3$) of the degree of $G_{n,p}$. 
Although our trees are somewhat smaller (by constant factor), we can handle trees
with degrees proportional to the minimum degree of the random graph.

\subsection{The algorithm}

All our results are proved using variants of the following very simple 
randomized embedding algorithm. First, choose arbitrarily some vertex $r$ 
of $T$ to be the {\em root}.
Then for every other vertex $u \in V(T)$ there is a unique path in $T$
from $r$ to $u$. The neighbor of $u$ on this path is called the {\em parent}
of $u$ and all the remaining neighbors of $u$ are called {\em children} of $u$.
The algorithms proceeds as follows. 

\paragraph{Algorithm 1.}
\noindent {\em Start by embedding the root $r$
at an arbitrary vertex $f(r) \in V(G)$. As long as $T$ is not completely
embedded, take an arbitrary vertex $u \in V(T)$ which is already embedded
but its children are not.
If $f(u)$ has enough neighbors in $G$ unoccupied by other vertices of $T$,
embed the children of $u$ by choosing vertices uniformly at random
from the available neighbors of $f(u)$ and continue. Otherwise, fail.}

\

This algorithm can be seen as a variant of a {\em tree-indexed random walk}, 
i.e. a random process corresponding to a tree where
each vertex assumes a random state depending only on the state of its parent.
The notion of a tree-indexed random walk was first introduced and studied
by Benjamini and Peres \cite{BP94}.
It is also used in the above mentioned paper of Benjamini and Schramm \cite{BS97}
to embed trees with a positive Cheeger constant into infinite expanding graphs.
In our case, we consider in fact a {\em self-avoiding}
tree-indexed random walk, where each state is chosen randomly,
conditioned on being distinct from previously chosen states.
The corresponding concept for a random walk is a 
well studied subject in probability (see, e.g., \cite{MS93}). Loosely speaking, 
we prove that our self-avoiding tree-indexed random walk behaves sufficiently randomly,
in the sense that it does not intersect the neighborhood of any vertex
more often than expected. To analyze the number of times the random process
intersects a given neighborhood, we use large deviation inequalities for supermartingales.

\subsection{A supermartingale tail estimate}
In all our proofs, we use the following tail estimate.

\begin{proposition}
\label{thm:supermartingale}
Let $X_1, X_2, \ldots, X_n$ be random variables in $[0,1]$ such that
for each $k$,
$$ \E[X_k \mid X_1, X_2,\ldots, X_{k-1}] \leq a_k.$$
Let $\mu = \sum_{i=1}^{n} a_i$. Then for any $0<\delta \leq 1$,
$$ \Pr[\sum_{i=1}^{n} X_i > (1+\delta) \mu] \leq
 e^{-\frac{\delta^2 \mu}{3}}.$$
\end{proposition}

This can be derived easily from the proof of Theorem 3.12(b) in~\cite{HMRR}.
We re-state this theorem here: {\em
Let $Y_1,Y_2,\ldots,Y_n$ be a martingale difference sequence with
$-a_k \leq Y_k \leq 1-a_k$ for each $k$, for suitable constants $a_k$;
and let $a = \frac{1}{n} \sum a_k$. Then for any $\delta>0$,
$$ \Pr[\sum_{k=1}^{n} Y_k \geq \delta a n] \leq
 e^{-\frac{\delta^2 a n}{2(1+\delta/3)}}.$$
}

A martingale difference sequence satisfies $\E[Y_i \mid Y_1,Y_2,\ldots,Y_{i-1}] = 0$.
However, it can be seen easily from the proof in \cite{HMRR} that
for this one-sided tail estimate, it is sufficient to assume 
$\E[Y_i \mid Y_1,Y_2,\ldots,Y_{i-1}] \leq 0$.
(Such a random process is known as a {\em supermartingale}.)
To show Proposition~\ref{thm:supermartingale}, set $Y_k = X_k - a_k$
 and $\mu = an = \sum_{i=1}^{n} a_k$.
The conditional expectations of $X_k$ are bounded by $a_k$, hence
the conditional expectations of $Y_k$ are non-positive as required.
Since $\delta \leq 1$, we also replace $2(1 + \delta/3)$ by $3$,
and Proposition~\ref{thm:supermartingale} follows.

Note also that we can always replace $\mu$ by a larger value
(e.g., by adding auxiliary random variables that are constants
with probability 1), and the conclusion still holds. Hence,
in Proposition \ref{thm:supermartingale}
it is enough to assume $\sum_{i=1}^{n} a_i \leq \mu$.

\section{Embedding trees in $C_4$-free graphs}
\label{section:C4-free}
The purpose of this section is to illustrate on a simple example
the main ideas and techniques that we will use
in our proofs. We start with $C_4$-free graphs,
which is a special case of two classes of graphs we are interested in:
graphs without short cycles, and graphs without $K_{s,t}$
(note that  $K_{2,2} = C_4$). 

Let's recall Algorithm 1. For a given rooted tree $T$,
we start by embedding the root $r \in V(T)$
at an arbitrary vertex $f(r) \in V(G)$. As long as $T$ is not completely
embedded, we take an arbitrary $u \in V(T)$ which is already embedded
but its children are not. If $f(u)$ has enough unoccupied neighbors in $G$,
we embed the children of $u$ uniformly at random in the available neighbors
of $f(u)$ and continue. Otherwise, we fail.

\begin{theorem}
\label{thm:C4-free}
Let $\epsilon \leq 1/8$, and let $G$ be $C_4$-free graph $G$
of minimum degree at least $d$. For any tree $T$ of size $|T| \leq \epsilon d^2$
and maximum degree $\Delta \leq d - 2 \epsilon d - 2$,
Algorithm 1 finds an embedding of $T$ in $G$ with high probability
(i.e., with probability tending to 1 when $d \rightarrow \infty$).
\end{theorem}

\paragraph{Example.} 
Before we plunge into the proof, let us consider the statement of this theorem
in a particular case, where $G$ is the incidence graph of a finite projective
plane. Let $q=d-1$ be a prime or a prime power and consider a 3-dimensional
vector space over the finite field $\mathbb{F}_q$. Let $V_1$ be all 2-dimensional linear subspaces of $\mathbb{F}^3_q$ (lines in a projective plane),
$V_2$ all 1-dimensional linear subspaces (points in a projective plane) 
and two vertices from $V_1$ and $V_2$ are adjacent if 
their corresponding subspaces contain one another.
This $G$ has $n = 2(q^2+q+1)=2(d^2 - d + 1)$ vertices, it is bipartite
and $d$-regular. Also, it is easy to see from the definition that $G$ contains
no $C_4$. Clearly, we cannot embed in $G$ trees of size larger than $O(d^2)$
or maximum degree larger than $d$. In this respect,
our theorem is tight up to constant factors.

It is also worth mentioning that in the analysis of our simple algorithm,
the trade-off between
the size of $T$ and the maximum degree $\Delta$ is close to being tight. Indeed, 
we show that for $\Delta = (1-\epsilon) d$, our algorithm cannot embed trees of size
much larger than $\epsilon d^2$. Suppose we are embedding a tree $T$
of depth $3$, where the degrees of the root and its children
(level $1$) are $\sqrt{d}$. On level $2$, the degrees are
$\epsilon d$ except one special vertex $z$ of degree $(1-\epsilon)d$.
On level $3$, there are only leaves. The size of this tree is
$\epsilon d^2 + \Theta(d)$.
 
We can assume that the root is embedded at a vertex corresponding to a point $a$.
The level-$1$ vertices are embedded into a set $L_1$ of $\sqrt{d}$ random lines
through $a$. The level-$2$ vertices are embedded into a set $P_2$ of
$d$ random points on these lines. Every point in the projective plane
(except $a$) has the same probability of appearing in $P_2$,
hence this probability is $d / (d^2 - d) = 1 / (d-1)$.
The level-$3$ vertices are embedded into random lines $L_3$ through
points in $P_2$, each line through a point in $P_2$ with probability $\epsilon$.
Now every line has probability roughly $\epsilon$ of being in $L_3$,
because one of its points on the average appears in $P_2$.
Consider the point where we embed the special vertex $z$ and assume
this is the last vertex we process in the algorithm. Each of the $d$ lines
through this point has probability roughly $\epsilon$ of being occupied
by a level-$3$ vertex, so on the average, only $(1-\epsilon)d$ lines
are available to host the children of $z$. Therefore, our algorithm
cannot succeed in embedding more than $(1-\epsilon) d$ children of $z$.

\paragraph{Proof of Theorem~\ref{thm:C4-free}.}
Let's fix an ordering in which the algorithm processes the vertices of $T$:
$V(T) = \{1,2,\ldots,|V(T)|\}$. Here, $1$ denotes the root and the ordering
is consistent with the structure of the tree in the sense that every vertex
can appear only after its parent. In step $0$, the algorithm embeds the root.
In step $t$, the children of $t$ are embedded randomly in the yet unoccupied
neighbors of $f(t) \in V(G)$. If $t$ is a leaf in $T$, the algorithm
is idle in step $t$.

Our goal is to argue that for large $d$, with high probability, 
the algorithm never fails.
The only way the algorithm can fail is that for a vertex $t \in V(T)$, embedded at 
$v = f(t) \in V(G)$, we are not able to place its children since too many neighbors
of $v$ in $G$ have been occupied by other vertices of $T$. This is the crucial
``bad event" we have to analyze:

{\em Let ${\cal B}_v$ denote the event that at some point, more than
$2 \epsilon d + 2$ neighbors of $v$ are occupied by vertices of $T$
other than the children of $f^{-1}(v)$.}

If we can show that with high probability, ${\cal B}_v$ does not occur 
for any $v \in V(G)$, then the algorithm clearly 
succeeds. To do this, we will modify our algorithm slightly and 
force it to stop immediately at the moment when the first bad event occurs. 
Thus, in analyzing ${\cal B}_v$, we can assume that for any $w \neq v$ the event
${\cal B}_w$ has not happened yet.
 
Our strategy is to prove that the probability of ${\cal B}_v$
for any given vertex $v$, even conditioned on our
embedding getting ``dangerously close" to $v$, is exponentially small in $d$.
Then, we argue that the number of vertices which can ever get dangerously
close to our embedding (i.e., the number of bad events we have to worry about)
is only polynomial in $d$. Therefore, we conclude that with high 
probability, no bad event
occurs.

\begin{lemma}
\label{lemma:bad-event}
Let $\epsilon \leq \frac18$ and $d \geq 24$. 
For a vertex $v \in V(G)$, condition on any history $\cal H$ of running
the algorithm up to a certain point such that at most $2$ vertices of $T$
have been embedded in $N(v)$. Then
$$ \Pr[{\cal B}_v \mid {\cal H}] \leq e^{-\epsilon d / 18}.$$
\end{lemma}

\noindent
{\bf Proof.}\,
For $t = 1,2,\ldots,|V(T)|$, let $X_t$ be an indicator variable of
the event that $f(t) \neq v$ but some child of $t$ gets embedded in $N(v)$.
Here we use the property that $G$ is $C_4$-free. Note that,
if $f(t) = w \neq v$, $w$ can have at most one neighbor in $N(v)$,
otherwise we get a $4$-cycle. Therefore, $t$ can have
at most one child embedded in $N(v)$ and $X_t$ represents the
number of vertices in $N(v)$, occupied by the children of $t$.

We condition on a history $\cal H$ of running the algorithm up to step $h$,
such that at most $2$ vertices of $N(v)$ have been occupied so far.
The bad event ${\cal B}_v$ can occur only if $X = \sum_{t=h+1}^{|T|} X_t 
> 2 \epsilon d$.  Therefore, our goal is to prove that 
this happens only with very small probability.

Each vertex chooses the embedding of its children randomly,
out of at least $d - 2 \epsilon d - 2$ still available choices
(here we assume that no bad event $B_w$ occurred before ${\cal B}_v$
for any $w \neq v$, or else the algorithm has failed already).
Thus we get
$$ \E[X_t] \leq \frac{d_T(t)}{d - 2 \epsilon d - 2} \leq 
 \frac{d_T(t)}{2d/3} $$
where $d_T(t)$ is the number of children of the vertex $t$ in $T$.
We also used $\epsilon \leq 1/8$ and $d \geq 24$.
This holds even conditioned on any previous history of the algorithm,
since the decisions for each vertex are made independently.
We are interested in the probability that $X = \sum_{t=h+1}^{|T|} X_t$
exceeds $2 \epsilon d$. 
Using the fact that $\sum_{t \in T} d_T(t) =|T|-1 \leq \epsilon d^2$,
we can bound the expectation of $X$ by
$$\mu = \E[X] = \sum_{t=h+1}^{|T|}\E[X_t] \leq
 \sum_{t \in T} \frac{d_T(t)}{2d/3} \leq \frac32 \epsilon d.$$ 
We use the supermartingale tail estimate
 (Proposition~\ref{thm:supermartingale}) with $\delta = \frac13$
 and $\mu = \frac32 \epsilon d$:
$$ \Pr[X > 2 \epsilon d] \leq e^{-\delta^2 \mu / 3}
 = e^{-\mu / 27} = e^{-\epsilon d / 18}.$$
Therefore, the bad event ${\cal B}_v$ happens with probability
at most $e^{-\epsilon d / 18}$. \hfill $\Box$

\vspace{2mm}

Our final goal is to argue that with high probability, no bad event ${\cal B}_v$
occurs for any vertex $v \in V(G)$. Since the number of vertices could
be potentially unbounded by any function of $d$, we cannot apply a 
straightforward union bound over all vertices in the graph. However,
we observe that the number
of vertices for which ${\cal B}_v$ can potentially occur is not very large.

Define ${\cal D}_v$ to be the event that at some point in the algorithm,
two vertices in $N(v)$ are occupied by vertices of $T$. This is the
event that the embedding of $T$ gets ``dangerously close" to $v$. Observe that if
${\cal D}_v$ is ``witnessed" by the pair of vertices of $T$ which are placed
in $N(v)$, each pair of vertices of $T$ can witness at most one event
${\cal D}_v$ (otherwise the same pair is in the neighborhood of two vertices
which implies a $C_4$). Since $T$ has at most $\epsilon d^2$ vertices,
the event ${\cal D}_v$ can occur for at most $\epsilon^2 d^4$ vertices
in any given run of the algorithm.

Clearly, event ${\cal B}_v \subseteq {\cal D}_v$.
Let's analyze the probability of ${\cal B}_v$, conditioned on ${\cal D}_v$. 
The event ${\cal D}_v$ can be written as a union of all histories $\cal H$
of running the algorithm up to the point where two vertices of $T$
get embedded in $N(v)$. By Lemma~\ref{lemma:bad-event},
$$ \Pr[{\cal B}_v \mid {\cal H}] < e^{-\epsilon d / 18} $$
for any such history $\cal H$. By taking the union of all these histories,
we get
$$ \Pr[{\cal B}_v \mid {\cal D}_v] < e^{-\epsilon d / 18}.$$
Now we can estimate the probability that ${\cal B}_v$ ever occurs for any vertex $v$:
\begin{eqnarray*}
\Pr[\exists v \in V; {\cal B}_v~\mbox{occurs}] & \leq & \sum_{v \in V} \Pr[{\cal B}_v]
 = \sum_{v \in V} \Pr[{\cal B}_v \mid {\cal D}_v] \Pr[{\cal D}_v]
 \leq e^{-\epsilon d / 18} \sum_{v \in V} \Pr[{\cal D}_v].
\end{eqnarray*}
Since ${\cal D}_v$ can occur for at most $\epsilon^2 d^4$ vertices in any given run
of the algorithm, we have $\sum_{v \in V} \Pr[{\cal D}_v] \leq \epsilon^2 d^4$. Thus
$$ \Pr[\exists v \in V; {\cal B}_v~\mbox{occurs}] \leq
 \epsilon^2 d^4 e^{-\epsilon d / 18} \rightarrow 0,$$
when $d \rightarrow \infty$.
Hence the algorithm succeeds with high probability.
\hfill $\Box$

\section{Embedding trees in $K_{s,t}$-free graphs}
\label{section:Kst-free}

Next, we consider the case of graphs which contain no complete bipartite subgraph $K_{s,t}$
with parts of size $s$ and $t$. We assume that $s \geq t$. It is known that the extremal size of such graphs depends
essentially only on the value of the smaller parameter $t$. Indeed, by the result of 
K\"ovari, S\'os and Tur\'an \cite{KST} the number of vertices 
in $K_{s,t}$-free graph with minimum degree $d$ is at least $c\,d^{t/(t-1)}$, where only
the constant $c$ depends on $s$. For relatively high values of
$s$ ($s >(t-1)!$) there are known constructions (see, e.g., \cite{KRS, ARS}) of
$K_{s,t}$-free graphs achieving this 
bound. Moreover, it is conjectured that $\Theta(d^{t/(t-1)})$ is the correct bound
for all $s \geq t$. This implies that one cannot 
embed trees larger than $O(d^{t/(t-1)})$ in a $K_{s,t}$-free graph
with minimum degree $d$. 
Also, it is obvious that the maximum degree in the tree should be $O(d)$. 
In this section we show how to embed trees with parameters very close to
these natural bounds that cannot be exceeded.
It is easier to analyze our algorithms in the case when the
maximum degrees in the tree are in fact bounded by $O(d/t)$.
First, we obtain this weaker result,
and then present a more involved analysis which shows that our algorithm
also works for trees with maximum degree at most $\frac{1}{256} d$.
Our algorithm here is a slight modification of Algorithm 1.

\paragraph{Algorithm 2.}
{\em For each vertex $v \in V(G)$, fix a set of $d$ neighbors
$N_+(v) \subseteq N(v)$. Start by embedding the root of the tree $r \in T$
at an arbitrary vertex $f(r) \in V(G)$. As long as $T$ is not
completely embedded, take an arbitrary vertex $u \in V(T)$ which is
already embedded but its children are not.
If $f(u)$ has enough neighbors in $N_+(f(u))$ unoccupied by other vertices of $T$,
embed the children of $u$ one by one, by choosing vertices uniformly at random
from the available vertices in $N_+(f(u))$, and continue. Otherwise, fail.}

\

The only difference from the original algorithm is that when embedding the children
of a vertex, we choose from a predetermined set of $d$ neighbors rather than
all possible neighbors. Since the maximum degree of $G$ can be very large, this 
modification is useful in the analysis of our algorithm. It allows us to bound
the number of dangerous events. However, we believe that the original
algorithm works as well and only our proof requires this modification.

\begin{theorem}
\label{thm:Kst-free}
Let $G$ be a $K_{s,t}$-free graph ($s \geq t$) with minimum degree $d$.
For any tree $T$ of size $|T| \leq \frac{1}{64} s^{-1/(t-1)} d^{t/(t-1)}$
and maximum degree $\Delta \leq \frac{1}{64 t} d$, Algorithm 2 finds
an embedding of $T$ in $G$ with high probability.
\end{theorem}

\noindent
{\bf Proof.}\,
We follow the strategy of defining bad events for each vertex $v \in V(G)$ and
bounding the probability that any such event occurs. 

{\em Let ${\cal B}_v$ denote the event that at some stage of the algorithm, more than
$\frac12 d + 2t$ vertices in $N_+(v)$ are occupied by vertices of $T$
other than children of $f^{-1}(v)$.}

Note that (as in the previous section), to bound the probability
of a bad event, we assume that our algorithm 
stops immediately at the moment when the first such event occurs.
To simplify our analysis, we also assume that the children of every vertex of $T$
are embedded in some particular order, one by one.
As long as ${\cal B}_v$ does not occur, we have at least $\frac12  d - 2t$
unoccupied vertices in $N_+(v)$.
Since degrees in the tree are bounded by $\frac{1}{64 t} d \leq \frac{1}{64} d$,
we have enough space for the children of any vertex 
to be embedded at $N_+(v)$. As we embed the children one by one,
the last child still has at least $\frac12 d - 2t - \frac{1}{64} d \geq \frac{1}{4} d$ 
choices available (for large enough $d$).

The new complication here is that another vertex $w$ could share many neighbors
with $v$. Unlike in the case of $K_{2,2}$-free graphs, where any two vertices can share
at most 1 neighbor, in $K_{s,t}$-free graphs (for $s > t \geq 2$),
we do not have any bound on the number of shared neighbors.
Therefore we have to proceed more carefully. For every vertex $v$ in $G$,
we partition all other vertices into two sets depending on how many neighbors they have in $N_+(v)$:
\begin{itemize}
\item $ L_v = \{w \neq v: |N_+(v) \cap N_+(w)| \leq 2 s^\frac{1}{t-1} d^\frac{t-2}{t-1} \}. $
\item $ M_v = \{w \neq v: |N_+(v) \cap N_+(w)| > 2 s^\frac{1}{t-1} d^\frac{t-2}{t-1} \}. $
\end{itemize}
The idea is that vertices in $L_v$ are harmless because the fraction of their children
that affects $N_+(v)$ is $O(d^{-1/(t-1)})$. Since the trees we are embedding
have size $O(d^{1 + 1/(t-1)})$, we show that the expected impact of these children on $N_+(v)$
is $O(d)$.

The vertices in $M_v$ have to be treated in a different way, because the 
fraction of their children in $N_+(v)$ could be very large. However, we 
prove that the total number of edges between $M_v$ and $N_+(v)$ cannot be 
too large, otherwise we would get a copy of $K_{s,t}$ in $G$. Therefore, 
the impact of the children of $M_v$ on $N_+(v)$ can be also controlled. 
Again, we ``start watching" a bad event for vertex $v$ only at the moment 
when it becomes dangerous.

{\em Let ${\cal D}_v$ denote the event that at least $t$ vertices in 
$N_+(v)$ are occupied by vertices of tree $T$ other than children of $f^{-1}(v)$.}

\begin{lemma}
\label{lemma:low-deg}
Let $\cal H$ be a fixed history of running the algorithm up to a point 
where at most $t$ vertices in $N_+(v)$ are occupied. 
Conditioned on $\cal H$, the probability that children of vertices embedded in $L_v$,
will ever occupy more than $\frac14 d+t$ vertices
in $N_+(v)$ is at most $e^{-d/24}$. 
\end{lemma}

\noindent
{\bf Proof.}\,
We use an argument similar to the proof of Lemma~\ref{lemma:bad-event}.
Fix an ordering of the vertices of $T$ starting from the root,
$i=1,2,\ldots,|T|$, as they are processed by the algorithm. 
Suppose that vertices $1, \ldots, h$ were embedded during the history $\cal H$. Let $X_i$ be the indicator variable of the event
that $i \in T$ is embedded in $N_+(v)$ and the parent of $i$ was embedded in $L_v$.
As long as the algorithm does not fail (i.e., no bad event happened), for each vertex $i \in T$ when it is 
embedded we have at least $d - \frac12 d - 2t - \frac{1}{64} d \geq \frac14 d$ choices
where to place the vertex. This holds even if we condition on any fixed embedding of 
vertices $j<i$. Moreover, the embedding decisions for different vertices are done independently.
Since we assume that the parent of $i$ was embedded in $L_v$,
at most $2 s^{1/(t-1)} d^{(t-2)/(t-1)}$ of these choices are in $N_+(v)$. Therefore,
conditioned on any previous history $\cal H$ such that $i$ was not embedded yet
$$ \Pr[X_i = 1 \mid {\cal H}] \leq \frac{2 s^\frac{1}{t-1}
 d^{\frac{t-2}{t-1}}}{\frac14 d} = 8 \left( \frac{s}{d} \right)^\frac{1}{t-1}.$$
Summing up over such vertices $i$ in the tree, whose number is at most
$|T| \leq \frac{1}{64} s^{-1/(t-1)} d^{t/(t-1)}$,
we have
$$ \E\Big[\sum_{i=h+1}^{|T|} X_i \mid {\cal H}\Big]
=\sum_{i=h+1}^{|T|} \E[X_i = 1 \mid {\cal H}] \leq |T| \cdot 8 \left( \frac{s}{d}
 \right)^\frac{1}{t-1} \leq \frac18 d.$$
Since, the upper bound on $\Pr[X_i = 1\mid {\cal H}]$ is still valid even if we also condition
on a fixed embedding of all vertices $j<i$, by Proposition~\ref{thm:supermartingale} with
$\mu = \frac18 d$ and $\delta = 1$,
$$ \Pr\Big[\sum_{i=h+1}^{|T|} X_i > \frac14 d \mid {\cal H}\Big] < e^{-d/24}. $$
By definition of $\cal H$, during the first $h$ steps of the algorithm
only at most $t$ vertices in $N_+(v)$ have been occupied.
Therefore, the probability that more than 
$\frac14 d + t$ vertices are ever occupied is at most $e^{-d/24}$. \hfill $\Box$

\medskip

Next, we treat the vertices whose parent is embedded in $M_v$. Recall that 
each vertex in $M_v$ has many neighbors in $N_+(v)$. However, the 
number of edges between $M_v$ and $N_+(v)$ cannot be too large. Observe 
that there is no $K_{s,t-1}$ in $G$ with $s$ vertices in $N_+(v)$ and 
$t-1$ vertices in $M_v$, otherwise we would obtain a copy of $K_{s,t}$ by adding $v$ 
to the part of size $t-1$. Also, this shows that for $t=2$, $M_v$ must be empty.
Indeed, by definition any vertex in $M_v$ has at least $2s$  neighbors in $N_+(v)$, which
together with vertex $v$ would form $K_{2s,2}$.
So in the following, we can assume $s \geq t \geq 3$.
The following is a standard estimate in extremal graph 
theory, whose short proof we include here for the sake of completeness.

\begin{lemma}
\label{lemma:Kst-extremal}
Consider a subgraph $H_v$ containing the edges between $M_v$ and $N_+(v)$, 
where $|N_+(v)| = d$, every vertex in $M_v$ has at least $2 s^{1/(t-1)} 
d^{(t-2)/(t-1)}$ neighbors in $N_+(v)$ and the graph does not contain 
$K_{s,t-1}$ (with $s$ vertices in $N_+(v)$ and $t-1$ vertices in $M_v$). 
Then $H_v$ has at most $2td$ edges.
\end{lemma}

\noindent
{\bf Proof.}\,
Let $m$ denote the number of edges in $H_v$ and assume $m > 2td$. Let $N$ 
denote the number of copies of $K_{1,t-1}$ (a star with $t-1$ edges) in $H_v$, with $1$ vertex in 
$N_+(v)$ and $t-1$ vertices in $M_v$. By convexity, the minimum number of $K_{1,t-1}$ in $H_v$ is 
attained when all vertices  in $N_+(v)$ have the same degree $m/|N_+(v)|$. Therefore
 $$ N \geq  |N_+(v)| {\frac{m}{|N_+(v)|} \choose t-1} = d {\frac{m}{d} \choose t-1}.$$
Our assumption that $m > 2td$ implies that $\frac{m}{d}, \frac{m}{d}-1,
\ldots, \frac{m}{d}-(t-2) \geq \frac{m}{2d}$ and therefore 
$$ N \geq d \frac{(\frac{m}{2d})^{t-1}}{(t-1)!}
 = \frac{m^{t-1}}{(t-1)! 2^{t-1} d^{t-2}}.$$ 
Since all the degrees in $M_v$ are at least $2 s^{1/(t-1)} d^{(t-2)/(t-1)}$, we have $m \geq 2 
s^{1/(t-1)} d^{(t-2)/(t-1)} |M_v|$. Then $m^{t-1} \geq 2^{t-1} s d^{t-2} 
|M_v|^{t-1}$ and
$$ N \geq \frac{m^{t-1}}{(t-1)! 2^{t-1} d^{t-2}} \geq \frac{s |M_v|^{t-1}}{(t-1)!} \geq s {|M_v| 
\choose t-1}.$$ 
Consequently, there must be a $(t-1)$-tuple in $M_v$ which 
appears in at least $s$ copies of $K_{1,t-1}$. This creates a copy of 
$K_{s,t-1}$, a contradiction.
 \hfill $\Box$

\begin{lemma}
\label{lemma:high-deg}
Let $\cal H$ be a fixed history of running the algorithm up to a point
where at most $t$ vertices in $N_+(v)$ are occupied. Then,
conditioned on $\cal H$, the probability that children of vertices 
embedded in $M_v$ will ever occupy more than $\frac14 d+t$ vertices
in $N_+(v)$ is at most $t \sqrt{d} e^{-\frac{1}{24t} \sqrt{d}}$.
\end{lemma}

\noindent
{\bf Proof.}\,
As we mentioned, we can assume $s \geq t \geq 3$, otherwise $M_v$ is empty.
Consider the vertices in $M_v$ and for every $w \in M_v$
denote the number of edges from $w$ to $N_+(v)$ by $d_w$.
We know that each vertex $w \in M_v$ has $d_w \geq 2 s^{1/(t-1)} d^{(t-2)/(t-1)} \geq 2 \sqrt{d}$
(using $t \geq 3$).  From Lemma~\ref{lemma:Kst-extremal},
we know that the total number of these edges is $\sum_{w \in M_v} d_w \leq 2td$.
This implies that $|M_v| \leq 2td / (2 \sqrt{d}) \leq t \sqrt{d}$.

For $w \in M_v$, let $X_w$ denote the number of tree vertices embedded in $N_+(v)$ after the
history $\cal H$, whose parent is embedded at $w$.
We claim that with high probability, $X_w \leq \frac{1}{8t} d_w$.
This can be seen as follows. Suppose that $f(x) = w$ for some $x \in V(T)$. The degree of $x$
in $T$ is at most $\frac{1}{64 t} d$ and the children of $x$ are embedded one by one.
Hence as we already explained, if no bad event ${\cal B}_w$ happened so far,  each child $y$ has at least $\frac 14 d$ choices
available for its embedding. Therefore, even conditioned on the embedding of the previous children,
the probability that $y$ is embedded in $N_+(v)$ is at most $p = \min \{1, d_w / (\frac14 d)\}$.
So $X_w$ satisfies the conditions of Proposition~\ref{thm:supermartingale} with
$\mu = \frac{1}{64 t} d \cdot d_w / (\frac14 d) = \frac{1}{16 t} d_w$.
By Proposition~\ref{thm:supermartingale} with $\delta = 1$,
$$ \Pr[X_w > \frac{1}{8t} d_w] \leq e^{-\mu/3} = e^{-\frac{1}{48t} d_w} \leq e^{-\frac{1}{24 t} \sqrt{d}}, $$
using $d_w \geq 2 \sqrt{d}$. By the union bound, the probability that $X_w > \frac{1}{8t} d_w$
for any $w \in M_v$ is at most $|M_v| e^{-\frac{1}{24t} \sqrt{d}} \leq t \sqrt{d} e^{-\frac{1}{24t} \sqrt{d}}$.
Otherwise, 
$$ \sum_{w \in M_v} X_w \leq \frac{1}{8t} \sum_{w \in M_v} d_w \leq \frac{1}{8t} \cdot 2td
 = \frac14 d.$$
Together with the $t$ vertices possibly occupied within history $\cal H$, this gives at most
$\frac14 d + t$ vertices occupied in $N_+(v)$.
\hfill $\Box$

Having finished all the necessary preparations we are now ready to 
complete the proof of Theorem~\ref{thm:Kst-free}. The bad event ${\cal B}_v$
can occur only if more than $\frac14 d + t$ vertices are occupied in 
$N_+(v)$ by children of vertices in $L_v$ or more than $\frac14 d + t$
vertices by children of vertices in $M_v$.
As we proved,  each of these events has probability smaller than $t \sqrt{d} e^{-\sqrt{d}/(24t)}$,
therefore the probability of ${\cal B}_v$ is at most $2t \sqrt{d} e^{-\sqrt{d}/(24t)}$.
This holds even if we condition on the event ${\cal D}_v$ (a disjoint union of histories $\cal H$)
which occurs  at the moment when $t$ vertices in $N_+(v)$ are occupied.

Let's estimate the number of events ${\cal D}_v$ which can occur. The event ${\cal D}_v$ is
witnessed by a $t$-tuple of vertices of tree $T$ which are embedded in $N_+(v)$.
The same $t$-tuple cannot be a witness to $s$ different events ${\cal D}_v$, because
then we would have a copy $K_{s,t}$ in our graph $G$. Therefore, each $t$-tuple can witness
at most $s-1$ events and the total number of events ${\cal D}_v$ is bounded by
$(s-1) |T|^t \leq s d^{2t}$. 
Since ${\cal D}_v$ can occur for at most $s d^{2t}$ vertices in any given run
of the algorithm, we have $\sum_{v \in V} \Pr[{\cal D}_v] \leq s d^{2t}$. 
Thus
\begin{eqnarray*}
\Pr[\exists v \in V; {\cal B}_v~\mbox{occurs}] & \leq & \sum_{v \in V} \Pr[{\cal B}_v]
 = \sum_{v \in V} \Pr[{\cal B}_v \mid {\cal D}_v] \Pr[{\cal D}_v] \\
 & \leq & 2t \sqrt{d} e^{-\frac{1}{24t} \sqrt{d}} \sum_{v \in V} \Pr[{\cal D}_v]
 \leq  2st \, d^{2t+\frac12} e^{-\frac{1}{24t} \sqrt{d}}
\end{eqnarray*}
which tends to $0$ as $d \rightarrow \infty$. \hfill $\Box$

\medskip

Finally, we show how to prove the same result for trees whose degrees can be
a constant fraction of $d$, independent of $t$. The following is a strengthened
version of Theorem~\ref{thm:Kst-free}.

\begin{theorem}
\label{thm:Kst-free2}
Let $G$ be $K_{s,t}$-free graph $G$ ($s \geq t$) of minimum degree $d$.
For any tree $T$ of size $|T| \leq \frac{1}{64} s^{-1/(t-1)} d^{t/(t-1)}$
and maximum degree $\Delta \leq \frac{1}{256} d$, Algorithm 2 finds
an embedding of $T$ in $G$ with high probability.
\end{theorem}

\noindent
{\bf Proof.}\,
The proof is very similar to the proof of Theorem~\ref{thm:Kst-free},
with some additional ingredients. We can assume that $t \geq 5$, otherwise
the result follows from Theorem~\ref{thm:Kst-free} directly.
We focus on the new issues arising from the
fact that degrees in the tree can exceed $O(d/t)$.
For a fixed vertex $v$, consider again the set $M_v$ defined by
 $$ M_v = \{w \neq v: |N_+(v) \cap N_+(w)| >
 2 s^\frac{1}{t-1} d^\frac{t-2}{t-1} \}. $$

We know from Lemma~\ref{lemma:Kst-extremal} that the number of edges from
$M_v$ to $N_+(v)$ is bounded by $2td$. Before, we argued that
since degrees are bounded by $O(d/t)$, the expected contribution of vertices
embedded along edges from $M_v$ to $N_+(v)$ cannot be too large.
The vertices in $T$ that could cause trouble are those embedded in $M_v$,
whose degree is more than $O(d/t)$. The contribution of the children
of these vertices to $N_+(v)$ might be too large. Hence we need to argue
that not too many vertices of this type can be embedded in $M_v$.

First, observe that using Lemma~\ref{lemma:Kst-extremal} and the
definition of $M_v$, the size of $M_v$ is bounded by
$$ |M_v| \leq \frac{e(M_v,N_+(v))}{2 s^{\frac{1}{t-1}} d^{\frac{t-2}{t-1}}}
\leq \frac{2td}{2 s^{\frac{1}{t-1}} d^{\frac{t-2}{t-1}}}
 \leq t d^{\frac{1}{t-1}}.$$
Similarly, if we denote by $Q$ the vertices of $T$ with degrees at least
$\frac{1}{64 t} d$, the number of such vertices is bounded by
$$ |Q| \leq \frac{2|T|}{\frac{1}{64 t} d}
 \leq \frac{\frac{1}{32} d^{\frac{t}{t-1}}}{\frac{1}{64 t} d}
 = 2 t d^{\frac{1}{t-1}}. $$
Our goal is to prove that not many vertices from $Q$ can be embedded in $M_v$.
For that purpose, we also need to define a new type of ``bad event" ${\cal C}_v$
and ``dangerous event" ${\cal E}_v$.

{\em The event ${\cal E}_v$ occurs if any vertex of the tree is embedded in $M_v$.
The event ${\cal C}_v$ occurs if after the first vertex embedded in $M_v$,
at least $8$ vertices from $Q$ are embedded in $M_v$.}

Now, consider any tree vertex $q \in Q$. At the moment when we embed $q$,
there are at least $\frac14 d$ choices, unless ${\cal B}_w$ happened for some 
vertex $w$ and the algorithm has failed already.
Since $|M_v| \leq t d^{\frac{1}{t-1}}$, the probability of embedding $q$
into $M_v$, even conditioned on any previous history $\cal H'$, is
$$ \Pr[f(q) \in M_v \mid {\cal H}'] \leq \frac{|M_v|}{\frac14 d} \leq
 \frac{4t d^{\frac{1}{t-1}}}{d} \leq \frac{4t}{d^{3/4}} $$
for $t \geq 5$.
We condition on any history $\cal H$ up to the first vertex embedded in $M_v$,
and estimate the probability that at least $8$ vertices from $Q$
are embedded in $M_v$ after this moment. For any particular $8$-tuple from $Q$,
this probability is bounded by $(4t / d^{3/4})^{8} = (4t)^8 / d^{6}$. 
The number of possible $8$-tuples in $Q$ is at most $|Q|^{8} \leq
 (2t d^{1/(t-1)})^{8} \leq (2t)^8 d^2$ for $t \geq 5$. Hence,
$$ \Pr[{\cal C}_v \mid {\cal H}] \leq \frac{(4t)^8}{d^{6}} (2t)^8 d^2
 = \frac{8^8 t^{16}}{d^4}.$$
By averaging over all histories up to the moment when the first vertex is
embedded in $M_v$, we get
$\Pr[{\cal C}_v \mid {\cal E}_v] \leq 8^8 t^{16} / d^4$.

Consider the number of events ${\cal E}_v$ that can ever happen.
For any event ${\cal E}_v$, there is a witness vertex $x \in V(T)$,
mapped to $f(x) = w \in M_v$. Observe that the definition of $w \in M_v$ is
symmetric with respect to $(v,w)$, i.e., we also have $v \in M_w$. We know that
$|M_w| \leq t d^{1 / (t-1)}$ for any $w \in V$, therefore each vertex
of the tree can be witness to at most $t d^{1 / (t-1)}$ events
${\cal E}_v$. In total, we can have at most $|T| \cdot t d^{1/(t-1)}
\leq d^{t/(t-1)} \cdot t d^{1/(t-1)} \leq t d^2$ events ${\cal E}_v$. 
Since ${\cal E}_v$ can occur for at most $t d^2$ vertices in any given run
of the algorithm, we have $\sum_{v \in V} \Pr[{\cal E}_v] \leq t d^2$. Hence,
\begin{eqnarray*}
\Pr[\exists v \in V; {\cal C}_v~\mbox{occurs}] & \leq & \sum_{v \in V} \Pr[{\cal C}_v]
 = \sum_{v \in V} \Pr[{\cal C}_v \mid {\cal E}_v] \Pr[{\cal E}_v] \\
 & \leq & \frac{8^8 t^{16}}{d^4} \sum_{v \in V} \Pr[{\cal E}_v] \leq 
 \frac{8^8 t^{16}}{d^4} t d^2 \leq \frac{8^8 t^{17}}{d^2}
\end{eqnarray*}
which tends to $0$ for $d \rightarrow \infty$. So, with high probability,
no event ${\cal C}_v$ happens.

Given that ${\cal C}_v$ does not occur for any vertex, we can carry out
the same analysis we used to prove Theorem~\ref{thm:Kst-free}.
The only difference is that each vertex $v$ might have up to $9$ vertices
from $Q$ embedded in $M_v$ ($8$ plus the first vertex ever embedded in $M_v$).
Since the degrees in $T$ are bounded by $\frac{1}{256} d$, even if 
the children of these vertices were embedded arbitrarily, still they can
occupy at most $\frac{9}{256} d$ vertices in $N_+(v)$. The number of vertices
in $N_+(v)$ occupied through vertices in $L_v$ or the contribution of the children
of vertices in $T$ with degree $O(d/t)$ that were embedded in $M_v$ can be analyzed just like in Theorem~\ref{thm:Kst-free}.
Thus, with high probability, at most $\frac12 d + \frac{9}{256} d + 2t < \frac34 d$
vertices are occupied in any neighborhood and so at least $\frac14 d$
vertices are always available to embed any vertex of the tree.  \hfill $\Box$

\section{Graphs of fixed girth}
\label{section:bounded-girth}

In this section we consider the problem of embedding trees into graphs
which have no cycle of length shorter than $2k+1$ for some $k > 1$.
(If the shortest cycle in a graph has length $2k+1$, such a graph
is said to have {\em girth} $2k+1$.) We also assume that the minimum degree
in our graph is at least $d$. It is easy to see that such $G$ must have
$\Omega(d^k)$ vertices,
because up to distance $k$ from any vertex $v$, $G$ looks locally like
a tree. It is widely believed that graphs of minimum degree $d$,
girth $2k+1$, and order $O(d^k)$ do exist for all fixed $k$ and large $d$.
Such constructions are known when $k = 2,3$ and $5$.
Since our graph might have order $O(d^k)$,  we cannot aspire to embed trees of size 
larger than $O(d^k)$ in $G$. This is what we achieve. 
For the purpose of analysis, we need to modify slightly our previous algorithms.

\paragraph{Algorithm 3.}
{\em For each $v \in V$, fix a set of its $d$ neighbors $N_+(v)$.
Assume that $T$ is a rooted tree with root $r$.
Start by making $k$ random moves from an arbitrary vertex $v_1 \in V$,
in each step choosing a random neighbor $v_{i+1} \in N_+(v_i)$.
Embed the root of the tree at $f(r) = v_k$.
 
As long as $T$ is not completely embedded, take an arbitrary vertex
$s \in V(T)$ which is embedded but its children are not.
If $f(s)$ has enough available neighbors in $N_+(f(s))$
unoccupied by other vertices of $T$, embed the children of $s$
among these vertices uniformly at random. Otherwise, fail.}

\medskip

The following is our main result for graphs of girth $2k+1$.

\begin{theorem}
\label{thm:fixed-girth}
Let $G$ be a graph of minimum degree $d$ and girth $2k+1$.
Then for any constant $\epsilon \leq \frac{1}{2k}$, Algorithm 3 succeeds
with high probability in embedding any tree $T$ of size
$ \frac14 \epsilon d^k$ and maximum degree $\Delta(T) \leq
 d - 2 \epsilon d - 2$.
\end{theorem}

To prove this theorem, we will generalize the analysis of the $C_4$-free
case to allow embedding of substantially larger trees. The solution is
to consider multiple levels of neighborhoods for each vertex. Starting
from any vertex $v \in V(G)$, we have the property that up to distance
$k$ from $v$, $G$ looks like a tree (otherwise we get a cycle of length
at most $2k$). Consequently, for any vertex $w$, there can be at most
one path of length $k$ from $w$ to $v$ . Therefore, embedding a subtree
whose root is placed at $w$ cannot impact the neighborhood of $v$ too much.

In fact, neighbors to be used in the embedding are chosen only
from a subset of $d$ neighbors $N_+(v)$. We can define an orientation
of $G$ where each vertex has out-degree exactly $d$, by orienting
all edges from $v$ to $N_+(v)$. (Some edges can be oriented both ways.)
Then, branches of the tree $T$ are embedded along {\em directed paths}
in $G$. 

\begin{definition}
For a rooted tree $T$, with a natural top-to-bottom orientation,
let $L_{k-1}(x)$ define the set of descendants $k-1$ levels down
from $x \in V(T)$.

For a tree vertex $x \in V(T)$, denote by $X_{v,x}$ the number of
vertices in $L_{k-1}(x)$ that end up embedded in $N_+(v)$,
before the children of $f^{-1}(v)$ are embedded.

For a vertex $v \in V(G)$, denote by $X_v$ the total number of vertices
in $T$ that end up embedded in $N_+(v)$,
before the children of $f^{-1}(v)$ are embedded.
\end{definition}

We extend $T$ to a larger rooted tree $T^*$ by adding a path of length $k-1$
above the root of $T$ and making the endpoint of this path
the root of $T^*$. Observe that our embedding algorithm proceeds effectively as if embedding $T^*$,
except the first $k-1$ steps do not occupy any vertices of $G$.
Each embedded vertex $y \in V(T)$ is a $(k-1)$-descendant of some
$x \in V(T^*)$ and hence $V(T) = \bigcup_{x \in V(T^*)} L_{k-1}(x)$.
By summing up the contributions over $x \in V(T^*)$, we get
$$ X_v = \sum_{x \in V(T^*)} X_{v,x}.$$

Our goal is to apply tail estimates on $X_v$ in order to bound
the probabilities of ``bad events". 
Just like before, we need to be careful in summing up these probabilities,
since the size of the graph might be too large for a union bound.
We start ``watching out" for the bad event ${\cal B}_v$ only after a
``dangerous event" ${\cal D}_v$ occurs. We also stop our algorithm immediately
after the first bad event happens.

{\em Event ${\cal B}_v$ occurs when $X_v > 2 \epsilon d + 2$.}
Event {\em ${\cal D}_v$ occurs whenever at least two vertices in $N_+(v)$
can be reached by directed paths of length at most $k-1$, avoiding $v$,
from the embedding of $T^*$. By the embedding of $T^*$, we also mean
the vertices visited in the first $k-1$ steps of the algorithm, which are not really occupied.}

Suppose $q_1,q_2$ are the first two vertices in $N_+(v)$ that can be reached by directed paths of length at 
most $k-1$, avoiding $v$, from the embedding of $T^*$.
Then we define a modified random variable $\tilde{X}_{v,x}$
as the number of vertices in $L_{k-1}(x)$, which are embedded in
$N_+(v) \setminus \{q_1,q_2\}$, but not through $v$ itself.
In other words, these random variables count the vertices
occupied in $N_+(v)$, not counting $q_1$ and $q_2$.
Observe that $X_v \leq \sum_{x \in V(T^*)} \tilde{X}_{v,x} + 2$.

\begin{lemma}
\label{lemma:girth-occupy}
Assume the girth of $G$ is at least $2k+1$.
Fix an ordering of the vertices of $T^*$ starting from the root,
$(x_1, x_2, x_3, \ldots)$, as they are processed by the algorithm.
Let $\cal H$ be a fixed history of running the algorithm until
two vertices $q_1, q_2 \in N_+(v)$ can be reached from an embedded vertex
by a directed path (avoiding $v$) of length at most $k-1$. Then for any vertex $x_i \in V(T^*)$,
$\tilde{X}_{v,x_i}$ is a $0/1$ random variable such that 
$$ \Pr[\tilde{X}_{v,x_i}=1 \mid {\cal H}, \tilde{X}_{v,x_1}, \tilde{X}_{v,x_2},
 \ldots, \tilde{X}_{v,x_{i-1}}]
 \leq \frac{|L_{k-1}(x_i)|}{(d-2\epsilon d - 2)^{k-1}}. $$
\end{lemma}

\noindent
{\bf Proof.}\,
First, note that any vertex $x_i$ embedded during the history $\cal H$ 
has $\tilde{X}_{v,x_i}=0$. (Since the only vertices in $N_+(v)$ possibly
reachable within $k-1$ steps from $f(x_j)$ are $q_1$ and $q_2$.)
Therefore we can assume that the embedding of $x_i$ together with the embedding
of the subtree of its descendants in $T^*$ is still undecided at the
end of ${\cal H}$. 
Let $\cal K$ denote the event that $x_i$ is embedded so that there
is a directed path of length exactly $k-1$  from $f(x_i)$ to $N_+(v)$, 
which avoids $v$ and has endpoint in $N_+(v)$ other than $q_1, q_2$.
Observe that this is the only way $\tilde{X}_{v,x_i}$ could be non-zero. 
Indeed, if $\tilde{X}_{v,x_i}=1$, then there is a branch of tree $T^*$ of length $k-1$
from $x_i$ to some $y$ that was mapped to a path from $f(x_i)$ to $N_+(v)$ such
that the vertex next to last is not $v$. However, such a path
from $f(x_i)$ to $N_+(v)$, if it
exists, is unique. If we had two different paths like this,
we could extend them to two paths of length $k$ between $f(x_i)$ and $v$,
which contradicts the girth assumption. Note that $\cal K$ occurs only
if this unique path leads to a vertex of $N_+(v)$ other than $q_1$ or $q_2$.
Also, we have that at most one vertex
$y \in L_{k-1}(x_i)$ can be embedded in $N_+(v)$. The variable
$\tilde{X}_{v,x_i}$ is equal to $1$ when this happens for some
$y \in L_{k-1}(x_i)$, and $0$ otherwise.

We bound the probability that $\tilde{X}_{v,x_i} = 1$, conditioned on
 $({\cal H}, \tilde{X}_{v,x_1}, \ldots, \tilde{X}_{v,x_{i-1}})$. In fact, let's 
condition even more strongly on a fixed
embedding $\cal E$ of all vertices of $T$ except for the descendants of $x_i$.
We also assume that $\cal E$ satisfies $\cal K$, i.e. $f(x_i)$ is
at distance exactly $k-1$ from $N_+(v)$, since otherwise
$\tilde{X}_{v,x_i}=0$. 
We claim that any such embedding implies
the values of $\tilde{X}_{v,x_1}, \ldots, \tilde{X}_{v,x_{i-1}}$. 
For vertices $x_j$ such that $L_{k-1}(x_j)$ does not intersect the subtree
of $x_i$, this is clear because the embedding of these vertices is fixed.
However, even if $L_{k-1}(x_j)$ intersects the subtree of $x_i$,
$\tilde{X}_{v,x_j}$ is still determined, since none of these vertices
can be embedded into $N_+(v)$. Indeed, any descendant of $x_i$ which is
in $L_{k-1}(x_j)$ must be also in $L_{k'}(x_i)$ for some $k' < k-1$.
If the embedding of $L_{k'}(x_i)$ intersects $N_+(v)$, we obtain
that there are two paths from $f(x_i)$ to $v$,
one of length $k$ and another of length $k'+1<k$. Together they form a cycle
of length shorter than girth, a contradiction.

Now fix a vertex $y \in L_{k-1}(x_i)$.
Every vertex $x_j \in T^*$, when embedded, chooses randomly from one of the available neighbors
of the vertex of $G$, in which its parent has been embedded.
As long as no bad event happened so far
(otherwise the algorithm would have terminated),
there are at least $d-2\epsilon d-2$ candidates available for $f(x_j)$.
Therefore, each particular vertex has probability at most
$1/(d-2\epsilon d-2)$ of being chosen to be $f(x_j)$. The probability that
$f(y) \in N_+(v)$ is the probability that our embedding follows
a particular path of length $k-1$. By the above discussion, this probability 
is at most $1/(d-2\epsilon d-2)^{k-1}$.
(Note that by our conditioning, this path might be already blocked by
the placement of other vertices; in such a case, the probability
is actually $0$.) Using the union bound, we have
$$ \Pr[\tilde{X}_{v,x_i}=1 \mid {\cal E}]  \leq
 \frac{|L_{k-1}(x_i)|}{(d - 2\epsilon d - 2)^{k-1}}. $$
Since the right hand side of this inequality is a constant, 
independent of the embedding, we get the same bound conditioned on 
$({\cal H}, \tilde{X}_{v,x_1}, \ldots, \tilde{X}_{v,x_{i-1}}, \cal K)$
and hence also conditioned on $({\cal H}, \tilde{X}_{v,x_1}, \ldots,
 \tilde{X}_{v,x_{i-1}})$. 
\hfill $\Box$

Now we are ready to use our supermartingale tail estimate from
Proposition~\ref{thm:supermartingale} to bound the
probability of a bad event.

\begin{lemma}
Assume $\epsilon \leq \frac{1}{2k}$ and $|T| \leq \frac14 \epsilon d^k$.
For any vertex $v \in V(G)$, condition on the dangerous event ${\cal D}_v$.
Then for large enough $d$, the probability that the bad event ${\cal B}_v$
happens is
$$ \Pr[{\cal B}_v \mid {\cal D}_v] \leq e^{-\epsilon d / 3}.$$
\end{lemma}

\noindent
{\bf Proof.}\,
The bad event means that $X_v > 2 \epsilon d + 2$. As before, first we
condition on any history $\cal H$ up to the point
when ${\cal D}_v$ happens. At this point, two vertices
$q_1,q_2 \in N_+(v)$ are within distance $k-1$ of the embedding of $T^*$
constructed so far. We consider these two vertices
effectively occupied.  Our goal is to prove that
the number of additional occupied vertices in $N_+(v)$ is small, 
namely $\sum_{i=1}^{|T^*|} \tilde{X}_{v,x_i} \leq 2 \epsilon d$.

By Lemma~\ref{lemma:girth-occupy}, we know that
$$ \Pr[\tilde{X}_{v,x_i} = 1 \mid {\cal H}, \tilde{X}_{v,x_1},
 \ldots, \tilde{X}_{v,x_{i-1}}] \leq
\frac{|L_{k-1}(x_i)|}{(d - 2\epsilon d - 2)^{k-1}}.$$
Therefore the expectation of $\tilde{X}_{v}=\sum_{i=1}^{|T^*|} \tilde{X}_{v,x_i} $ is bounded by
$$ \mathbb{E}[\tilde{X}_{v}]=\sum_{i=1}^{|T^*|} \mathbb{E}[\tilde{X}_{v,x_i}]
 \leq \sum_{i=1}^{|T^*|} \frac{|L_{k-1}(x_i)|}{(d-2\epsilon d - 2)^{k-1}}
 \leq \frac{|T|}{(d-2\epsilon d - 2)^{k-1}} < \frac{4|T|}{d^{k-1}}
  \leq \epsilon d. $$
Here we used that $\epsilon \leq \frac{1}{2k}$, $d$ large enough, and
 $|T| \leq \frac14 \epsilon d^k$.
So we can set $\mu = \epsilon d$, $\delta = 1$ and use Proposition~\ref{thm:supermartingale} to conclude that,
$$ \Pr\big[ \tilde{X}_{v} > 2 \epsilon d \mid {\cal H}\big]
 \leq e^{-\epsilon d / 3}. $$
The same holds when we condition on the event ${\cal D}_v$,  which is the disjoint union
of all such histories $\cal H$. Consequently, 
$X_v \leq \tilde{X}_{v} + 2 \leq 2 \epsilon d + 2$
with high probability, which concludes the proof.
\hfill $\Box$

\medskip

To finish the proof of Theorem~\ref{thm:fixed-girth},
we show that with high probability, ${\cal B}_v$ does not happen 
for any vertex $v \in V$. First, let's examine how many events ${\cal D}_v$ can possibly occur for a given run
of the algorithm. Every vertex $v$ for which ${\cal D}_v$ happens has a ``witness
 pair" of vertices in $N_+(v)$ satisfying the condition that they can be reached
by directed paths of length at most $k-1$ from the embedding of $T^*$.
The number of such vertices is at most $|T^*| d^{k-1} \leq d^{2k}$.
Also, observe that the same pair can be a witness to at most
$1$ event ${\cal D}_v$, otherwise we have a 4-cycle in $G$ which 
contradicts the high girth property. 
Hence the number of possible witness pairs is at most 
$$ {d^{2k} \choose 2} \leq d^{4k} $$
and each event ${\cal D}_v$ has a unique witness pair.
Therefore, the expected number of events ${\cal D}_v$ is
$$ \sum_v \Pr[{\cal D}_v] \leq d^{4k}.$$
Now we bound the probability that any bad event ${\cal B}_v$ occurs.
\begin{eqnarray*}
\Pr[\exists v \in V; {\cal B}_v~\mbox{occurs}] & \leq & \sum_{v \in V} \Pr[{\cal B}_v]
 = \sum_{v \in V} \Pr[{\cal B}_v \mid {\cal D}_v] \Pr[{\cal D}_v] \\
 & \leq & e^{-\epsilon d / 3} \sum_{v \in V} \Pr[{\cal D}_v] \leq
  d^{4k} e^{-\epsilon d / 3}.
\end{eqnarray*}
For a constant $k$ and $d \rightarrow \infty$, this probability
tends to $0$. \hfill $\Box$

\section{Random graphs and the property ${\cal P}(d,k,t)$}
\label{section:bounded-paths}

The main objective of this section is to obtain nearly optimal tree embedding results
for random graphs. In our analysis, we do not actually require true randomness.
The important condition that $G$ has to satisfy is a certain ``pseudorandomness"
property, stated below.  Roughly speaking, the property requires
that there are not too many paths between any pair of vertices,
compared to how many paths a random graph would have.

\paragraph{Property ${\cal P}(d,k,t)$.}
Let $d, k$ and $t$ be positive integers. A graph $G$ on $n$ vertices satisfies property ${\cal P}(d,k,t)$ if
\begin{enumerate}
\item $G$ has minimum degree at least $d$.
\item For any $u,v \in V$, the number of paths of length $k$ from $u$ to $v$ is 
$$ P_{k}(u,v) \leq d^{1/4}.$$
\item For any $u,v \in V$, the number of paths of length $k+1$ from $u$ to $v$ is
 $$ P_{k+1}(u,v) \leq \frac{d^{k+1}}{t}.$$
\end{enumerate}

\paragraph{Remark.}
In the second condition, $d^{1/4}$ is somewhat arbitrary.
For $k$ constant, it would be enough to require $P_{k}(u,v) =
 o(d / \log d)$. However, having a larger gap between $P_{k}(u,v)$
and $d$ allows our framework to work for larger (non-constant)
values of $k$.

Observe that $d$-regular graphs of girth $2k+1$  satisfy ${\cal P}(d,k,t=d^k)$,
because there is at most one path of length $k$ between any pair of vertices.
Thus our embedding results for graphs satisfying this property  
implies similar statements for regular graphs of fixed girth,
although somewhat weaker than those we presented in Section~\ref{section:bounded-girth}.
Our main focus in this section is on random graphs.

\begin{proposition}
\label{random}
A random graph $G_{n,p}$ where $\frac12 \geq p \geq n^{a-1}$, $a > 0$ constant,
satisfies almost surely ${\cal P}(d,k,t)$ with $t = (1-o(1)) n$, $d=(1-o(1))pn$
and $k \geq 1$ chosen so that
$$\frac{1}{4} (pn)^{-3/4} < p^k n^{k-1} \leq \frac{1}{4} (pn)^{1/4}.$$
\end{proposition}

\noindent
{\bf Proof.}\,
Since we assume $pn \geq n^a$, we have $k \leq 1 + 1/a$,
otherwise $p^k n^{k-1} = p (pn)^{k-1} \geq pn >> (pn)^{1/4}$
contradicting our choice of $k$. Hence, $k$ is a constant.

The degree of every vertex in $G_{n,p}$ is a binomially distributed
random variable with parameters $n$ and $p$. Thus, by standard tail
estimates (Chernoff bounds), the probability that it is smaller than
$$ d = pn - \sqrt{pn} \log n = (1-o(1)) pn $$
is $e^{-\Omega(\log^2 n)} = o(1/n)$. Therefore with high probability
the minimum degree of $G_{n,p}$ is at least $d$.

The expected number of paths of length $k$ from $u$ to $v$ is
$$ \E[P_{k}(u,v)] \leq p^k n^{k-1} \leq \frac{1}{4} (pn)^{1/4} $$
by our choice of $k$. We use the Kim-Vu inequality \cite{KimVu}
to argue that $P_k(u,v)$ is strongly concentrated. Let $t_e$
be the indicator variable of edge $e$. We can write
$$ P_k(u,v) = \sum_{P} \prod_{e \in P} t_e $$
where $P$ runs over all possible paths of length $k$ between $u$ and $v$.
Clearly, this is a multilinear polynomial of degree $k$.
Let $\frac{\partial}{\partial t_I} P_k(u,v)$ denote the partial derivative of $P_k(u,v)$ with respect to all variables in the 
set $I$. Using the notation of \cite{KimVu}, we set
$$ E_i = \max_{|I|=i} \E\left[\frac{\partial}{\partial t_I} P_k(u,v)\right],$$
$E = \max_{i \geq 0} E_i$ and $E' = \max_{i \geq 1} E_i$. In particular, 
$E_0$ is the expected value of $P_k(u,v)$. 
The Kim-Vu inequality states that
$$ \Pr\big[|P_k(u,v) - E_0| > a_k \lambda^k \sqrt{E'E}\big]
 = O\big(e^{-\lambda + (k-1) \log n}\big) $$
for any $\lambda > 1$ and $a_k = 8^k \sqrt{k!}$.
In our case, $\E\left[\frac{\partial}{\partial t_I} P_k(u,v)\right]$
can be seen as the expected number of $u$-$v$ paths of length $k$ with $i$ edges
already fixed to be on the path. For any choice of such $i$ edges,
if $i<k$, we have at most $n^{k-i-1}$ choices to complete the path and
the probability that such a path appears is $p^{k-i}$. Hence,
$E_i \leq  p^{k-i} n^{k-i-1}$ for $i<k$. For $i=k$, we have $E_k = 1$.
Hence, $E = \max_{i \geq 0} E_i \leq p^k n^{k-1} \leq \frac14 (pn)^{1/4}$ and
$E' = \max_{i \geq 1} E_i \leq 1$.  By the Kim-Vu inequality with
$\lambda = (k+2) \log n$, we have
$$ \Pr\big[|P_k(u,v) - E_0| > a'_k (pn)^{1/8} \log^k n\big] =
 O\big(e^{-3 \log n}\big) = O\big(n^{-3}\big),$$
where $a'_k=(k+2)^k a_k=8^k (k+2)^k \sqrt{k!}$.
Thus, we get for all pairs $(u,v)$ that with high probability
$$P_k(u,v) \leq E_0 + a'_k (pn)^{1/8} \log^k n
 \leq \frac14 (pn)^{1/4} +a'_k (pn)^{1/8} \log^k n <\frac12 (pn)^{1/4} \leq d^{1/4}.$$

To estimate $P_{k+1}(u,v)$, we use a similar argument. Again, this is
a multilinear polynomial $P_{k+1}(u,v) = \sum_P \prod_{e \in P} t_e$,
this time of degree $k+1$. The expectation is
$E_0 = \E[P_{k+1}(u,v)] \leq p^{k+1} n^k$. 
Further, we get $E_i \leq p^{k+1-i} n^{k-i}$ for $i < k$,
$E_{k+1} = 1$ and therefore, $E = \max_{i \geq 0} E_i = E_0$.
Since our choice of $k$ implies that $E_0=(1-o(1))p^{k+1} n^k > (pn)^{1/4}/5$, we
also have  
$$E' = \max_{i \geq 1} E_i = \max\big(p^{k} n^{k-1}, 1\big) \leq 5E_0/(pn)^{1/4}.$$
By Kim-Vu with $\lambda = (k+2) \log n$,
$$ \Pr\big[|P_{k+1}(u,v) - E_0| > a'_k \sqrt{EE'} \log^k n\big] =
 O\big(e^{-3 \log n}\big) = O\big(n^{-3}\big),$$
where $a'_k =(k+2)^ka_k$ is a constant. Note that 
$a'_k \sqrt{EE'} \log^k n \leq 5a'_k\log^k n E_0/(pn)^{1/8}=o(E_0)$. Recall also that 
$d=(1-o(1))pn$ and $t=(1-o(1))n$. Thus, for all pairs $(u,v)$ with high probability
$$\hspace{3.5cm} P_{k+1}(u,v) \leq E_0 + o(E_0)
 \leq (1+o(1)) p^{k+1} n^k \leq d^{k+1} / t \hspace{3.5cm} \Box$$

\paragraph{Algorithm 4.}
{\em 
Start by making $k$ random moves from an arbitrary vertex $v_0 \in V$,
in each step choosing a random neighbor $v_{i+1} \in N(v_i)$.
Embed the root of the tree $r \in T$ at $f(r) = v_k$.

As long as $T$ is not completely embedded, take an arbitrary vertex
$u \in V(T)$ which is embedded but its children are not.
If $f(u)$ has enough available neighbors in $N(f(u))$
unoccupied by other vertices of $T$, embed the children of $u$
one by one by choosing vertices randomly from the available neighbors
of $f(u)$. Otherwise, fail.}

\medskip
The following is our main theorem.

\begin{theorem}
\label{thm:bounded-paths}
Let $G$ be a graph on $n$ vertices satisfying property ${\cal P}(d,k,t)$ 
for $d \geq \log^8 n$, $k \leq \log n$ and $\epsilon, \delta > 0$ are
such that 
\begin{equation}
\label{eq:eps-delta}
(2k \epsilon)^{1/k} + \delta + \frac{1}{k} \leq 1.
\end{equation}
Then for any tree $T$ of maximum degree at most $\delta d$ and size
at most $\epsilon t$, the algorithm above 
finds embedding of $T$ with high probability.
\end{theorem}

\noindent
This result has an interesting consequence already for $k=1$.
Let $G$ be a graph on $n$ vertices with 
minimum degree $pn$ such that every two distinct vertices of $G$ have
at most $O(p^2 n)$ common neighbors. For $p \gg n^{-1/2}$ there are several known explicit construction of such graphs 
and their properties were extensively studied by various researchers
(see, e.g., survey \cite{KS1} and its references).
Our theorem implies nearly optimal embedding results for such $G$ and shows that it contains every tree of 
order $\Omega(n)$ with maximum degree $\Omega(pn)$.

Considering the extreme values of $\epsilon$ and $\delta$
that satisfy (\ref{eq:eps-delta}), we obtain embeddings of
\begin{itemize}
\item trees with maximum degree at most a constant fraction of $d$
 (e.g., $\frac{1}{4} d$) and size $2^{-\Theta(k)} t$.
\item trees with maximum degree $O(d/k)$ and size $O(t/k)$.
\end{itemize}
Combining Theorem \ref{thm:bounded-paths} with Proposition \ref{random},
we see that for a random graph $G_{n,p}$ with $p = n^{a-1}$ and
constant $a > 0$  we can use $d \simeq pn$, $t \simeq n$
and $k \simeq  1/a$. Therefore for such $p$ we are embedding trees whose size and maximum degree 
are proportional to the order and minimum degree of $G_{n,p}$.
This is clearly tight up to constant factors.

\medskip

Before proving the theorem, we outline the strategy of our proof.
Our goal is to argue that there is some $\alpha > 0$ such that
no more than $\alpha d$ vertices are ever occupied in any neighborhood $N(v)$,
including vertices embedded through $v$ itself.
Again, we consider the number $X_v$ of vertices in $N(v)$ occupied
by vertices of $T$, other than those embedded as children of $v$.
The ``bad event" ${\cal B}_v$ occurs when $X_v > d/k$ and we stop the algorithm immediately after the 
first such event.
At most $\delta d$ vertices can be embedded as children of $v$,
therefore assuming that no bad event happens, 
at most $(1/k + \delta) d$ vertices are eventually occupied in any neighborhood
$N(v)$.  Since $1/k + \delta \leq 1 - (2 k \epsilon)^{1/k}$ by (\ref{eq:eps-delta}),
we can set
$$\alpha = 1 - (2 k \epsilon)^{1/k}.$$
If no bad even occurs, any vertex of $T$ has at least $(1-\alpha) d$ choices available
for its embedding. If a bad event occurs, we can assume that the algorithm fails.

We estimate the probability of ${\cal B}_v$ by studying the random variable $X_v$.
The expectation $\E[X_v]$ is bounded relatively easily,
since this is determined by the number of possible ways
that a vertex of $T$ can reach the neighborhood $N(v)$.
This can be bounded using our property ${\cal P}(d,k,t)$.
The more challenging part of the proof is to argue that the probability of ${\cal B}_v$
is very small, since the contributions from different vertices of the tree
are not independent. We handle this issue by dividing the contributions
into blocks of variables which are effectively independent.
We write $X_v = \sum_{i=1}^{k} Y_{v,i}$ and use a supermartingale
tail estimate to bound each $Y_{v,i}$.

The following definitions are similar to those in Section~\ref{section:bounded-girth}.

\begin{definition}
For a rooted tree $T$, with a natural top-to-bottom orientation,
let $L_{k-1}(x)$ define the set of descendants $k-1$ levels down
from $x \in V(T)$.

For a vertex $v \in V(G)$, denote by $X_v$ the number of vertices
in $T$ that end up embedded in $N(v)$, before the children of $f^{-1}(v)$
are embedded. 
 
For a tree vertex $x \in V(T)$, denote by $X_{v,x}$ the number of
vertices in $L_{k-1}(x)$ that end up embedded in $N(v)$,
before the children of $f^{-1}(v)$ are embedded.
\end{definition}

As in Section~\ref{section:bounded-girth}, we extend $T$ to a larger tree $T^*$
by adding a path of $k$ auxiliary vertices above the root. 
Each embedded vertex $y$ is a $(k-1)$-descendant of some
$x \in V(T^*)$ and hence $V(T) = \bigcup_{x \in V(T^*)} L_{k-1}(x)$.
By summing up the contributions over $x \in V(T^*)$, we get
$$ X_v = \sum_{x \in V(T^*)} X_{v,x}.$$

\begin{lemma}
\label{lemma:1-subtree}
Assume $G$ satisfies property ${\cal P}(d,k,t)$ and fix a tree vertex $x \in V(T)$.
Then $X_{v,x}$ is bounded by $d^{1/4}$ with probability $1$, and
$$ \E[X_{v,x} \mid {\cal T}] \leq (1-\alpha)^{-k}  |L_{k-1}(x)| \frac{d}{t} $$
where $\cal T$ is any fixed embedding of the entire tree $T$ except
for the vertex $x$ and its descendants.
\end{lemma}

\noindent
{\bf Proof.}\,
Assume that conditioned on $\cal T$, the parent $q$ of $x$
is embedded at $f(q) = w \in V(G)$. The only way that a vertex $y \in L_{k-1}(x)$
can end up in $N(v)$ (but not through $v$) is when some branch of the tree $T$ from $q$ to $y$ is
embedded in a path of length $k$ from $w$ to $N(v)$, avoiding $v$.
Such paths can be extended uniquely to paths of length $k+1$ from $w$
to $v$.  We know that the number of such paths
is bounded by $P_{k+1}(w,v) \leq d^{k+1} / t$.

Since there are at least $(1-\alpha) d$ choices when we embed
each vertex, the probability of following a particular path
of length $k$ is at most $\frac{1}{((1-\alpha)d)^{k}}$.
By the union bound, the probability that $y$ is embedded in $N(v)$
is 
$$ \Pr[f(y) \in N(v) \mid {\cal T}] \leq
  \frac{1}{((1-\alpha)d)^{k}} P_{k+1}(w,v)
 \leq \frac{d}{(1-\alpha)^{k} t}.$$ Finally,
$$ \E[X_{v,x} \mid {\cal T}] = \sum_{y \in L_{k-1}(x)}
 \Pr[f(y) \in N(v) \mid {\cal T}] \leq
  \frac{|L_{k-1}(x)| d}{(1-\alpha)^{k} t}.$$

Similarly, the number of paths of length $k-1$ from any vertex $u$ to $N(v)$, avoiding $v$,
is the same as the number $P_k(u,v)$ of paths of length $k$ from $u$ to $v$.
Even if all these $P_{k}(u,v)$ paths are used in the embedding of $T$,
the vertices in $L_{k-1}(x)$ cannot occupy more than $P_{k}(u,v)$ neighbors of $v$.
Therefore, we can always bound $X_{v,x} \leq P_{k}(u,v) \leq d^{1/4}$.
\hfill $\Box$

Next, we want to argue about the concentration of $X_v = \sum_{x \in V(T^*)}
X_{v,x}$. Since the placements of different vertices in $T$ are highly
correlated, it is not clear whether any concentration result applies
directly to this sum. However, we can circumvent this obstacle by
partitioning $V(T^*)$ into subsets where the dependencies can work
only in our favor.

\begin{definition}
Let $r^*$ be the root of $T^*$, then every vertex of $T^*$ is in
$L_j(r^*)$ for some $j$. Define a partition $V(T^*) = W_0 \cup W_1 \cup
 \ldots \cup W_{k-1}$ by
$$ W_j = \bigcup_{j'=j\pmod k} L_{j'}(r^*).$$
For each vertex $v \in V(G)$ and $0 \leq j < k$, define
$$ Y_{v,j} = \sum_{x \in W_j} X_{v,x}.$$
\end{definition}

Obviously, we have $X_v = \sum_{x \in V(T^*)} X_{v,x} = \sum_{j=0}^{k-1} Y_{v,j}$.
In the following, we argue that each $Y_{v,j}$ has a very small one-sided tail.

\begin{lemma}
\label{lemma:one-class}
Let $\ell_j = \sum_{x \in W_j} |L_{k-1}(x)|$. Then 
 $\E[Y_{v,j}] \leq (1-\alpha)^{-k} \frac{\ell_j}{t} d$ and
$$ \Pr\left[Y_{v,j} > (1-\alpha)^{-k}
 \left(  \frac{\ell_j}{t} + \frac{\epsilon}{k} \right) d \right]
 < e^{-\frac{\epsilon d^{3/4}}{3k^2 (1-\alpha)^{k}}}.$$
\end{lemma}

\noindent
{\bf Proof.}\,
By Lemma~\ref{lemma:1-subtree}, we know that
$ \E[X_{v,x} \mid {\cal T}] \leq (1-\alpha)^{-k}  |L_{k-1}(x)| \frac{d}{t}$
where $\cal T$ is  any fixed embedding of $T$ except $x$ and its
subtree. Therefore, the same also holds without any conditioning.
By taking a sum over all $x \in W_j$,
$$ \E[Y_{v,j}] = \sum_{x \in W_j} \E[X_{v,x}] \leq (1-\alpha)^{-k} 
 \sum_{x \in W_j} |L_{k-1}(x)| \frac{d}{t} =
  (1-\alpha)^{-k} \ell_j \frac{d}{t}.$$
For a tail estimate, we use Proposition~\ref{thm:supermartingale}.
Write the vertices of $W_j= \{x_1, x_2, \ldots, x_r \}$ in order as they are embedded by the algorithm
and write $X_i = d^{-1/4} X_{v,x_i}$.
The important observation is that the values of $X_1, X_2, \ldots, X_{i-1}$
are determined if we are given the embedding of the tree $T$ except for the vertex $x_i$ and its subtree
(let's denote this condition by ${\cal T}_i$). This holds because
$X_1, \ldots, X_{i-1}$ depend only on the embedding of vertices
$x_1, \ldots, x_{i-1}$ and their subtrees of depth $k-1$.
Since all these vertices are either at least $k$ levels above $x_i$
in the tree $T$, or on the same level or below (but not in the subtree of $x_i$),
their subtrees of depth $k-1$ are disjoint from the subtree of $x_i$.
Hence, conditioning on ${\cal T}_i$ is stronger than conditioning on
$X_1,\ldots,X_{i-1}$. Since $\E[X_i \mid {\cal T}_i] = d^{-1/4} \E[X_{v,x_i} \mid {\cal T}_i]
 \leq  (1-\alpha)^{-k} |L_{k-1}(x_i)| \frac{d^{3/4}}{t}$, we can also write
$$ \E[X_i \mid X_1,\ldots,X_{i-1}] \leq  (1-\alpha)^{-k} |L_{k-1}(x_i)| \frac{d^{3/4}}{t}. $$
The range of $X_{v,x_i}$ is $[0, d^{1/4}]$, hence $X_i \in [0,1]$.
Summing over $W_j$, we have  $\sum \E[X_i] \leq (1-\alpha)^{-k} \ell_j
 \frac{d^{3/4}}{t}$, so let's set $\mu = (1-\alpha)^{-k} \ell_j \frac{d^{3/4}}{t}$.
By Proposition~\ref{thm:supermartingale},
$$ \Pr[\sum X_i > (1 + \epsilon') \mu]  < e^{-\frac{\epsilon'^2 \mu}{3}}.$$
Using that $\ell_j \leq |T| \leq \epsilon t$, for $\epsilon' = \frac{t \epsilon}{\ell_j k}$, we get
$$ \Pr\left[\sum X_i > \mu +  \frac{\epsilon}{k}(1-\alpha)^{-k} d^{3/4} \right]
 < e^{-\frac{t \epsilon^2 d^{3/4}}{3 \ell_j k^2 (1-\alpha)^{k}}}
 < e^{-\frac{\epsilon d^{3/4}}{3k^2 (1-\alpha)^k}}. $$
Since $Y_{v,j} = d^{1/4} \sum X_i$, this proves the claim of the lemma.
\hfill $\Box$

\begin{lemma}
\label{lemma:gen-bad-event}
Let ${\cal B}_v$ denote the ``bad event" that $X_v > d/k$.
Assuming that (\ref{eq:eps-delta}) holds,
and $|T| \leq \epsilon t$, then for any fixed vertex $v \in V$
the bad event happens with probability
$$ \Pr[{\cal B}_v] < k e^{-\frac{d^{3/4}}{6 k^3}}.$$
\end{lemma}

\noindent
{\bf Proof.}\,
We have $X_v = \sum_{j=0}^{k-1} Y_{v,j}$.
Recall that $(1-\alpha)^k = 2k \epsilon$.
By Lemma~\ref{lemma:one-class},
$$ \Pr\left[Y_{v,j} > (1-\alpha)^{-k} \left( \frac{\ell_j}{t}  + \frac{\epsilon}{k} \right) d \right]
 <   e^{-\frac{\epsilon d^{3/4}}{3 k^2 (1-\alpha)^{k}}}  = e^{-\frac{d^{3/4}}{6 k^3}} $$
for each $j=0,1,2,\ldots,k-1$.
By the union bound, the probability that any of these events happens is
at most $k e^{-{d^{3/4}}/{6 k^3}}$. If none of them happen,  we have
$$ X_v = \sum_{j=0}^{k-1} Y_{v,j} \leq 
 (1-\alpha)^{-k} \sum_{j=0}^{k-1} \left( \frac{\ell_j}{t}
   + \frac{\epsilon}{k} \right) d
 = (1-\alpha)^{-k} \left( \frac{|T|}{t} + \epsilon \right) d
 \leq (1-\alpha)^{-k} \cdot 2 \epsilon d = \frac{d}{k}.$$
\hfill $\Box$

To finish the proof of Theorem~\ref{thm:bounded-paths}, we note that
$d \geq \log^8 n$ and $k \leq \log n$. The probabilities of bad events
${\cal B}_v$ are bounded by $k e^{-d^{3/4} / 6 k^3} \leq (\log n) \ e^{-\frac16 \log^3 n}
 \leq 1/n^{\log n}$.
There are $n$ potential bad events, so none of them occurs with high
probability.

\section{Concluding remarks}

In this paper we have shown that a very simple randomized algorithm
can find efficiently tree embeddings with near-optimal parameters, surpassing
some previous results achieved by more involved approaches.
Here are few natural questions which remain open. 

\begin{itemize}
\item
It would be interesting to extend our results from graphs of girth $2k+1$ to
graphs without cycles of length $2k$. For $k=3$, this follows
from our work combined with a result of Gy\"{o}ri. In \cite{G97}
he proved that every bipartite
$C_6$-free graph can be made also $C_4$-free by deleting
at most half of its edges. Therefore given a $C_6$-free graph
with minimum degree $d$, we can first take its maximum bipartite
subgraph. This will decrease the number of edges by at most factor of two.
Then we can use the above mentioned result of Gy\"{o}ri to obtain
a $C_4$-free and $C_6$-free graph which has at least a quarter of the
original edges, i.e., average degree at least $d/4$.
In this graph we can find a subgraph  where the minimum degree is at least $d/8$ 
($1/2$ of average degree). Since it is bipartite, this subgraph has no cycles of length shorter than $7$.
This shows that every $C_6$-free graph $G$ with minimum degree $d$ contains
a subgraph $G'$ of girth at least $7$ whose minimum degree is a constant fraction
of $d$. Using our result, we can embed in $G'$ (and hence also in $G$)
every tree of size $O(d^3)$ and maximum degree $O(d)$.

More generally, it is proved in \cite{KO05} that any $C_{2k}$-free
graph contains a $C_4$-free subgraph with at least $\frac{1}{2(k-1)}$-fraction
of its original edges. Moreover it is conjectured in \cite{KO05},
that any $C_{2k}$-free graph contains a subgraph of girth $2k+1$
with at least an $\epsilon_k$-fraction of the edges.
If this conjecture is true, it shows that the tree embedding
problems for $C_{2k}$-free graphs and graphs of girth $2k+1$ are
equivalent up to constant factors.

\item
For random graphs $G_{n,p}$ our approach works most efficiently when
the edge probability $p=n^{a-1}$ for some constant $a>0$.
Nevertheless, it can be used to embed trees in sparser random graphs as well.
By analyzing more carefully the application of the Kim-Vu inequality,
one can show that for every fixed $\epsilon > 0$, a random graph with edge 
probability $p \geq e^{\log^{1/2+\epsilon} n}/n$ satisfies ${\cal P}(d,k,n/2)$
with $d \simeq pn$ and $k \simeq \log_d n$. However, when $p= n^{-1+o(1)}$ we 
have $k \rightarrow \infty$ and therefore both the maximum degree an the size of the tree we can embed are only an $o(1)$-fraction of the optimum. 
It would be extremely interesting to show that for edge probability
$p=n^{-1+o(1)}$, perhaps even $p=c/n$ for some large constant $c>0$,
the random graph $G_{n,p}$  still contains every tree with maximum degree
$O(pn)$ and size $O(n)$.

It would be also nice to weaken our pseudorandomness property
${\cal P}(d,k,t)$ which is defined in terms of numbers of paths
between pairs of vertices. The most common definition of pseudorandomness
is in terms of edge density between subsets of vertices of a graph. 
In particular, it would be interesting to extend our results to embedding
of trees in graphs whose edge distribution is close to that of random graph. 

\item 
Finally, we wonder if there are any additional interesting families of graphs
for which one can show that our simple randomized algorithm succeeds
to embed trees with nearly optimal parameters. 
\end{itemize}


\end{document}